\documentclass{amsart}

\usepackage{amsmath,amssymb,amsthm}
\usepackage{graphicx,tikz}
\newtheorem{theorem}{Theorem}

\newtheorem*{proposition}{Proposition}

\newtheorem*{lemma}{Lemma}

\theoremstyle{definition}

\theoremstyle{remark}

\begin{document}

\title[]{Improved Bounds for Hermite-Hadamard inequalities in higher dimensions}

\author[]{Thomas Beck}
\address{Department of Mathematics, University of North Carolina at Chapel Hill, CB\#3250 Phillips Hall,
Chapel Hill, NC 27599}
\email{tdbeck@email.unc.edu}

\author[]{Barbara Brandolini}
\address{Universit\`a degli Studi di Napoli Federico II, Dipartimento di Matematica 
e Applicazioni ``R. Caccioppoli'', Complesso Monte S. Angelo - Via Cintia, 80126 Napoli, Italia.}
\email{brandolini@unina.it}

\author[]{Krzysztof Burdzy}
\address{Department of Mathematics, Box 354350, University of Washington, Seattle, WA
98195}
\email{ burdzy@uw.edu}

\author[]{Antoine Henrot}
\address{Institut Elie Cartan de Lorraine, CNRS UMR 7502 and Universite de 
Lorraine, BP 70239 54506 Vandoeuvre-les-Nancy, France}
\email{ antoine.henrot@univ-lorraine.fr}

\author[]{Jeffrey J. Langford}
\address{Bucknell University, 1 Dent Drive, Lewisburg, PA 17837, USA}
\email{Jeffrey.Langford@bucknell.edu}

\author[]{Simon Larson}
\address{Department of Mathematics, KTH Royal Institute of Technology, SE-100 44
Stockholm, Sweden}
\email{larson.simon@gmail.com}

\author[]{Robert Smits}
\address{Department of Mathematical Sciences, New Mexico State University,
Las Cruces, NM 88003-8001}
\email{rsmits@nmsu.edu}

\author[]{Stefan Steinerberger}
\address{Department of Mathematics, Yale University, New Haven, CT 06511, USA}
\email{stefan.steinerberger@yale.edu}

\begin{abstract} Let $\Omega \subset \mathbb{R}^n$ be a convex domain and let $f:\Omega \rightarrow \mathbb{R}$ be a positive, subharmonic function (i.e. $\Delta f \geq 0$). Then
$$ \frac{1}{|\Omega|} \int_{\Omega}{f dx} \leq \frac{c_n}{ |\partial \Omega| } \int_{\partial \Omega}{ f d\sigma},$$
where $c_n \leq 2n^{3/2}$. This inequality was previously only known for convex functions with a much larger constant. We also
show that the optimal constant satisfies $c_n \geq n-1$. As a byproduct,
we establish a sharp geometric inequality for two convex domains where one contains the other $ \Omega_2 \subset \Omega_1 \subset \mathbb{R}^n$:
 $$ \frac{|\partial \Omega_1|}{|\Omega_1|} \frac{| \Omega_2|}{|\partial \Omega_2|} \leq n.$$
\end{abstract}

\maketitle

\section{Introduction}
\subsection{Convex functions.} The Hermite-Hadamard inequality dates back to an 1883 observation of Hermite \cite{hermite} with an independent use by Hadamard \cite{hada} in 1893: it says that
for convex functions $f:[a,b] \rightarrow \mathbb{R}$ 
$$ \frac{1}{b-a} \int_{a}^{b}{f(x) dx} \leq \frac{f(a) + f(b)}{2}.$$
This inequality is rather elementary and has been refined in many ways -- we refer to the monograph of Dragomir \& Pearce \cite{drago}. However, there is relatively little work outside of the one-dimensional case; we refer to \cite{cal1, cal2, jianfeng, nicu2, nicu, choquet, past, stein1}.
 The strongest possible statement that one could hope for is, for convex functions $f:\Omega \rightarrow \mathbb{R}$ defined on convex domains $\Omega \subset \mathbb{R}^n$,
$$ \frac{1}{|\Omega|} \int_{\Omega}{f ~d x} \leq \frac{1}{|\partial \Omega|} \int_{\partial \Omega}{f ~d\sigma}.$$
This inequality has been shown to be true for many special cases: it is known for $\Omega = \mathbb{B}_3$ the $3-$dimensional ball by Dragomir \& Pearce \cite{drago} and $\Omega = \mathbb{B}_n$ by de la Cal \& Carcamo \cite{cal1} (other proofs are given by de la Cal, Carcamo \& Escauriaza \cite{cal2} and Pasteczka \cite{past}), the simplex \cite{bes}, the square \cite{square}, triangles \cite{chen} and Platonic solids \cite{past}. 
It was pointed out by Pasteczka \cite{past} that if the inequality holds for a domain $\Omega$ with constant 1, then plugging in affine functions shows that the center of mass of $\Omega$ and the center of mass of $\partial \Omega$ coincide, which is not generally true for convex bodies; therefore the inequality cannot hold with constant 1 in higher dimensions uniformly over all convex bodies. 
The first uniform estimate was shown in \cite{stein1}: if $f:\Omega \rightarrow \mathbb{R}$ is a convex, positive function on the convex domain $\Omega \subset \mathbb{R}^n$, then we have 
\begin{equation} \label{full}
 \frac{1}{|\Omega|} \int_{\Omega}{f ~dx} \leq  \frac{c_n}{|\partial \Omega|} \int_{\partial \Omega}{f~d\sigma}
 \end{equation}
 with $c_n \leq 2 n^{n+1}$.
In this paper, we will improve this uniform estimate and show that the optimal constant satisfies $n-1 \leq c_n \leq 2 n^{3/2}$. We do not have a characterization of the extremal convex functions $f$ on a given domain $\Omega$ (however, see below, we have such a characterization in the larger family of subharmonic functions).

\subsection{Subharmonic functions.}
Niculescu \& Persson \cite{choquet} (see also \cite{cal2, nicu2}) have pointed out that one could also seek such inequalities for subharmonic functions, i.e. functions satisfying $\Delta f \geq 0$. We note that all convex functions are subharmonic. Jianfeng Lu and the last author \cite{stein1} showed that for all positive, subharmonic functions $f:\Omega \rightarrow \mathbb{R}$ on convex domains $\Omega \subset \mathbb{R}^n$
\begin{equation} \label{jianfeng}
\int_{\Omega}{f ~d x} \leq |\Omega|^{1/n} \int_{\partial \Omega}{f ~d \sigma}.
\end{equation}
Estimates relating the integral of a positive subharmonic function $f$ over $\Omega$ to the integral over the boundary $\partial \Omega$ are linked to the torsion function on $\Omega$ given by
\begin{align*}
- \Delta u &= 1 \quad \mbox{in}~\Omega\\
u &= 0 \quad \mbox{on}~\partial \Omega.
\end{align*}
Integration by parts and the inequalities $u \geq 0, \Delta f \geq 0$ show that
\begin{align*}
\int_{\Omega}{ f dx} = \int_{\Omega}{ f (-\Delta u) dx} &= \int_{\partial \Omega}{ \frac{\partial u}{\partial \nu} f d\sigma} - \int_{\Omega}{(\Delta f) u dx} \\
&\leq  \int_{\partial \Omega}{ \frac{\partial u}{ \partial \nu} f d\sigma} \\
&\leq \max_{x \in \partial \Omega}{  \frac{\partial u}{\partial \nu}(x)} \int_{\partial \Omega} { f d\sigma},
\end{align*}
where $\nu$ is the inward pointing normal vector.
This computation suggests that we may have the following characterization of the optimal constant for a given convex domain $\Omega$.
\begin{proposition}[see e.g. \cite{dragomirk, choquet}] The optimal constant $c(\Omega)$ in the inequality 
$$\int_{\Omega}{f ~d x} \leq c(\Omega) \int_{\partial \Omega}{f ~d \sigma}$$
for positive subharmonic functions is given by 
$$ c(\Omega) = \max_{x \in \partial \Omega}{  \frac{\partial u}{\partial \nu}(x)}.$$
\end{proposition}
The lower bound on $c(\Omega)$ follows from setting $f$ to be the Poisson extension of a Dirac measure located at the point at which the normal derivative assumes its maximum. The derivation also shows that it suffices to consider the case of harmonic functions $f$. 
Implicitly, this also gives a characterization of extremizing functions (via the Green's function).
 Jianfeng Lu and the last author \cite{jianfeng} used this proposition in combination with a gradient estimate for the torsion function to show that the best constant in (\ref{jianfeng}) is uniformly bounded in the dimension.
We will follow a similar strategy to obtain an improved bound for the optimal constant in (\ref{jianfeng}).

\section{The Results}
Our first result improves the constant $c_n$ from (\ref{full}) in all dimensions for subharmonic functions and shows that the growth is at most polynomial.  
\begin{theorem} Let $\Omega \subset \mathbb{R}^n$ be convex and let $f:\Omega \rightarrow \mathbb{R}$ be a positive, subharmonic function. Then
\begin{equation} \label{thm1}
 \frac{1}{|\Omega|} \int_{\Omega}{f dx} \leq \frac{c_n}{ |\partial \Omega| } \int_{\partial \Omega}{ f d\sigma},
 \end{equation}
where the optimal constant $c_n$ satisfies
$$ c_n \leq \begin{cases} n^{3/2}~\,\,\,\,\,\,\,\,\,\,\,\,\,\,\,\, \mbox{if}~n~\mbox{is odd}, \\
\frac{n^2+n}{\sqrt{n+2}}  ~\qquad\mbox{if}~n~\mbox{is even}. \end{cases}$$
\end{theorem}

In particular, for $n=2$ dimensions, our proof shows the inequality
$$ \frac{1}{|\Omega|} \int_{\Omega}{f dx} \leq \frac{3}{ |\partial \Omega| } \int_{\partial \Omega}{ f d\sigma},$$
where the constant 3 improves on constant 8 obtained earlier for convex functions in \cite{stein1}. 
To complement the result in Theorem 1 we prove that any constant for which \eqref{thm1} is valid must grow at least linearly with the dimension. 
\begin{theorem} The optimal constant $c_n$ in (\ref{thm1}) is non-decreasing in $n$ and satisfies
\begin{equation} \label{convex}
 c_n \geq \max\{n-1, 1\}.
 \end{equation}
\end{theorem}
In order to prove Theorem 2 we establish a connection to an isoperimetric problem that is of interest in its own right. Specifically, we prove the following Lemma.
\begin{lemma} In any dimension $n\geq 1$,
\begin{equation} \label{convex}
 \sup \left\{ \frac{|\partial \Omega_1|}{|\Omega_1|} \frac{| \Omega_2|}{|\partial \Omega_2|}: \Omega_2 \subset \Omega_1 ~ \emph{both convex domains in}~\mathbb{R}^{n}\right\}=n.
 \end{equation}
\end{lemma}
We are not aware of any prior treatment of this shape optimization problem in the literature.
Problem (\ref{convex}) can be equivalently written as
$$ \sup \left\{ \frac{|\partial \Omega|}{|\Omega|} \frac{1}{h(\Omega)}: \Omega ~ \mbox{a convex set in}~\mathbb{R}^{n-1}\right\},~ \mbox{where}~ h(\Omega) = \inf_{X \subset \overline{\Omega}}{ \frac{|\partial X|}{|X|}}$$
 denotes the Cheeger constant. We refer to Alter \& Caselles \cite{alter} and Kawohl \& Lachand-Robert \cite{kawohl}. 
The result of Kawohl \& Lachand-Robert \cite{kawohl}
will be a crucial ingredient in the proof of Theorem 2 (we note that the infimum runs over all subsets, it is known that the Cheeger set is unique and convex).

We also obtain a slight improvement of the constant in (\ref{jianfeng}).

\begin{theorem} Let $\Omega \subset \mathbb{R}^n$ be convex and let $f:\Omega \rightarrow \mathbb{R}$ be a positive, subharmonic function. Then
$$ \int_{\Omega}{f dx} \leq \frac{ |\Omega|^{1/n}}{ \omega_n^{1/n} \sqrt{n} } \int_{\partial \Omega}{ f d\sigma},$$
where $\omega_n$ is the volume of the unit ball in $n-$dimensions.
\end{theorem}
We observe that, as $n$ tends to infinity, $\omega_n^{1/n} \sqrt{n} \rightarrow \sqrt{2 \pi e}.$
We also note a construction from \cite{jianfeng} which shows that 
the constant in Theorem 3 is at most a factor $\sqrt{2}$ from optimal in high dimensions.


\section{Proof of Theorem 1}
\subsection{Convex functions.}
We first give a proof of Theorem 1 under the assumption that $f$ is convex; this argument is fairly elementary and is perhaps useful in other settings. A full proof of Theorem 1 is given in \S 3.2.

\begin{proof}
This proof combines three different arguments. The first argument is that
\begin{equation} \label{b1}
  \int_{\Omega}{f dx} \leq \frac{w(\Omega)}{2} \int_{\partial \Omega}{ f d\sigma}
  \end{equation}
from the one-dimensional Hermite-Hadamard inequality applied along fibers that are orthogonal to the hyperplanes realizing the width $w(\Omega)$. 

\begin{figure}[h!]
\begin{center}
\begin{tikzpicture}[scale=1.6]
\draw [very thick] (0,0) -- (5,0);
\draw [dashed] (-0.5,0) -- (5.5,0);
\draw [dashed] (-0.5,1.05) -- (5.5,1.05);
\draw [very thick] (0,0) to[out=45, in=135] (5,0);
\draw [thick, <->] (6,0) -- (6,1);
\node at (6.3, 0.5) {$w(\Omega)$};
\node at (2.5, 0.5) {$\Omega$};
\draw [dashed] (1,0) -- (1,0.7);
\draw [dashed] (1.1,0) -- (1.1,0.7);
\draw [dashed] (1.2,0) -- (1.2,0.75);
\draw [dashed] (1.3,0) -- (1.3,0.8);
\end{tikzpicture}
\end{center}
\caption{Application of the one-dimensional inequality on a one-dimensional fiber. This step is lossy if the boundary is curved.}
\end{figure}
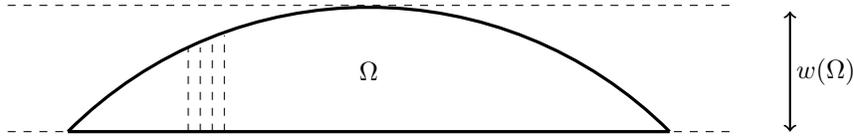

Steinhagen \cite{steinhagen} showed that width can be bounded in terms of the inradius
\begin{equation} \label{b2}
w(\Omega) \leq \begin{cases} 2 \sqrt{n} \cdot \mbox{inrad}(\Omega)~\,\,\,\,\,\qquad \mbox{if}~n~\mbox{is odd,} \\
2 \frac{n+1}{\sqrt{n+2}} \cdot\mbox{inrad}(\Omega) ~\qquad\mbox{if}~n~\mbox{is even}. \end{cases}
\end{equation}
The last inequality follows from \cite{larson1}: if $\Omega \subset \mathbb{R}^n$ is a convex body and
\begin{equation*}
	\Omega_t = \{x\in \Omega: d(x, \partial\Omega)>t\},
\end{equation*}
where $d(x, \partial \Omega)$ denotes the distance to the boundary
$$ d(x, \partial \Omega) = \inf_{y \in \partial \Omega}{ \|x-y\|},$$
then
$$|\partial \Omega_t| \geq |\partial \Omega| \left(1 - \frac{t}{\mbox{inrad}(\Omega)}\right)_+^{n-1}.$$
Since $|\nabla d(x, \partial \Omega)| = 1$ almost everywhere, the coarea formula implies
\begin{align*} |\Omega| &= \int_{0}^{ \tiny \mbox{inrad}(\Omega)} |\partial \Omega_t| dt \\
&\geq |\partial \Omega|\int_{0}^{ \tiny \mbox{inrad}(\Omega)}{ \left(1 - \frac{t}{\mbox{inrad}(\Omega)}\right)^{n-1} dt} = |\partial \Omega| \frac{\mbox{inrad}(\Omega)}{n}
\end{align*}
and thus we obtain (also stated in \cite[Eq. 13]{larson2})
\begin{equation} \label{b3}
\mbox{inrad}(\Omega) \leq n \frac{|\Omega|}{|\partial \Omega|}.
\end{equation}
Combining inequalities (\ref{b1}), (\ref{b2}) and (\ref{b3}) implies the result.
\end{proof}
Both Steinhagen's inequality as well as inequality (\ref{b3}) are sharp for the regular simplex. However, see Fig. 1, an application of the one-dimensional Hermite-Hadamard inequality can only be sharp if the fibers are hitting the boundary at a point at which they are normal, otherwise there is a Jacobian determinant determined by the slope of the boundary and better results are expected. It is not clear to us how to reconcile these two competing factors.

\subsection{A Proof of Theorem 1}
\begin{proof}  We have, for all positive, subharmonic functions $f:\Omega \rightarrow \mathbb{R}$,
$$ \int_{\Omega}{ f dx}  \leq \max_{x \in \partial \Omega}{  \frac{\partial u}{\partial \nu}(x)} \int_{\partial \Omega} { f d\sigma},$$
where $u$ is the torsion function.
A classic bound on the torsion functions is given in Sperb \cite[Eq. (6.12)]{sperb}),
$$  \max_{x \in \partial \Omega}{  \frac{\partial u}{\partial \nu}(x)} \leq \sqrt{2} \|u\|^{\frac{1}{2}}_{L^{\infty}}.$$
Moreover, using Steinhagen's inequality in combination with (\ref{b3}), we know that $\Omega$ is contained within a strip of width
$$ w(\Omega) \leq   \frac{|\Omega|}{|\partial \Omega|} \begin{cases} 2 n^{3/2} ~\,\,\,\,\,\qquad \mbox{if}~n~\mbox{is odd,} \\
2 \frac{n^2+n}{\sqrt{n+2}} ~\qquad\mbox{if}~n~\mbox{is even}. \end{cases}$$
We can now use the maximum principle to argue that the torsion function in $\Omega$ is bounded from above by the 
torsion function in the strip of width $w(\Omega)$ (see Fig. 2). That torsion function, however, is easy to compute since the problem
becomes one-dimensional. Orienting the strip to be given by
$$ S = \left\{(x,y) \in \mathbb{R}^{n-1} \times \mathbb{R}: |y| \leq \frac{w(\Omega)}{2} \right\},$$
we see that the torsion function on the strip is given by
$$ v(x,y) = \frac{w(\Omega)^2}{8} - \frac{y^2}{2}.$$

\begin{figure}[h!]
\begin{center}
\begin{tikzpicture}[scale=0.8]
\draw[very thick] (0,0) ellipse (3cm and 1cm);
\draw[thick] (-4, 1) -- (4, 1);
\draw[thick] (-4, -1) -- (4, -1);
\draw [<->] (5, -1) -- (5,1);
\node at (5.5,0) {$w(\Omega)$};
\node at (0,0) {$\Omega$};
\end{tikzpicture}
\end{center}
\caption{The torsion function in $\Omega$ is bounded from above by the torsion function of the strip.}
\end{figure}
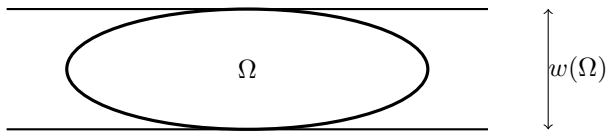

This shows
$$\|u\|^{}_{L^{\infty}} \leq \frac{w(\Omega)^2}{8}$$
and thus
\begin{align*}
  \max_{x \in \partial \Omega}{  \frac{\partial u}{\partial \nu}(x)}  &\leq \sqrt{2} \|u\|_{L^{\infty}}^{1/2} \\
 &\leq \frac{w(\Omega)}{2}  \leq  \frac{|\Omega|}{|\partial \Omega|} \begin{cases}  n^{3/2} ~\,\,\,\,\,\qquad \mbox{if}~n~\mbox{is odd} \\
 \frac{n^2+n}{\sqrt{n+2}} ~\qquad\mbox{if}~n~\mbox{is even}. \end{cases}
\end{align*}
\end{proof}

\section{Proof of Theorem 2}

The purpose of this section is to prove $c_{n+1} \geq c_n$ as well as the inequality
$$  c_n \geq \sup \left\{ \frac{|\partial \Omega_1|}{|\Omega_1|} \frac{| \Omega_2|}{|\partial \Omega_2|}: \Omega_2 \subset \Omega_1 ~ \mbox{both convex domains in}~\mathbb{R}^{n-1}\right\}.$$
Theorem 2 is then implied by this statement together the proof of the Geometric Lemma in Section \S 5.

\begin{proof}
The proof is based on explicit constructions. We first show that $c_{n+1} \geq c_{n}$. This is straightforward and based on an extension in the $(n+1)-$first coordinate: for any $\varepsilon>0$, we can find a convex domain $\Omega_{\varepsilon} \subset \mathbb{R}^n$
and a positive, convex function $f_{\varepsilon}:\Omega_{\varepsilon} \rightarrow \mathbb{R}$ such that
$$ \frac{1}{|\Omega_{\varepsilon}|} \int_{\Omega_{\varepsilon}}{f_{\varepsilon}  dx} \geq \frac{c_n - \varepsilon}{ |\partial \Omega_{\varepsilon}| } \int_{\partial \Omega_{\varepsilon}}{ f_{\varepsilon}  d\sigma}.$$
We define, for any $z > 0$, 
$$ \Omega_{z,\varepsilon} = \left\{(x, y): x \in \Omega_{\varepsilon} ~\mbox{and}~ 0 \leq y \leq z\right\} \subset \mathbb{R}^{n+1}$$
and $f_{z, \varepsilon}:\Omega_{z, \varepsilon} \rightarrow \mathbb{R}$ via 
$$ f_{z, \varepsilon}(x,y) = f_{\varepsilon}(x).$$
Then 
$$ \frac{1}{|\Omega_{z,\varepsilon}|} \int_{\Omega_{z,\varepsilon}}{f_{z,\varepsilon}  dx dy} =  \frac{1}{|\Omega_{\varepsilon}|} \int_{\Omega_{\varepsilon}}{f_{\varepsilon}  dx}.$$
This integral simplifies for $z$ large since
$$ \lim_{z \rightarrow \infty}  \frac{1}{ |\partial \Omega_{z,\varepsilon}| } \int_{\partial \Omega_{z,\varepsilon}}{ f_{z,\varepsilon}  d\sigma} =   \frac{1}{ |\partial \Omega_{\varepsilon}| } \int_{\partial \Omega_{\varepsilon}}{ f_{\varepsilon}  d\sigma}.$$
Picking $\varepsilon$ sufficiently small and $z$ sufficiently large shows that $c_{n+1} < c_n$ leads to a contradiction.\\

We now establish the inequality inequality 
$$  c_n \geq \sup \left\{ \frac{|\partial \Omega_1|}{|\Omega_1|} \frac{| \Omega_2|}{|\partial \Omega_2|}: \Omega_2 \subset \Omega_1 ~ \mbox{both convex domains in}~\mathbb{R}^{n-1}\right\}.$$
To this end pick $0 \in \Omega_2 \subset \Omega_1 \subset \mathbb{R}^{n-1}$ in such a way that both domains are convex. We will now define a domain $\Omega_N \subset \mathbb{R}^n$ and a convex function $f_N:\Omega_N \rightarrow \mathbb{R}$ where $N \gg 1$ will be a large parameter. We first define the convex sets
$$ C_1 = \left\{(x, y): x \in \Omega_1 \,~\mbox{and}\,~ y \geq - N^3\right\}$$
and
$$ C_2 = \left\{(x,y): x \in \left(1- \frac{y}{N^2}\right)\Omega_2 \,~\mbox{and}~  y \leq N \right\}.$$
The set $\Omega_N$ is then given as the intersection $\Omega_N = C_1 \cap C_2$ (see Fig. 3). We observe that $\Omega_N$ is the intersection of two convex sets and is therefore convex. Also, looking at the scaling, we see that $C_1$ dominates: looking at $\Omega_N$ from `far away`, it looks essentially like $C_1$ truncated. We now make this precise: note that there exists a constant $\lambda \geq 1$ such that $\Omega_1 \subseteq \lambda \Omega_2$ and then
$$ \Omega_N \cap \left\{(x,y) \in \mathbb{R}^{n}: y \leq -(\lambda-1) N^2 \right\} =  C_1 \cap \left\{(x,y) \in \mathbb{R}^{n}: y \leq -(\lambda-1) N^2 \right\}.$$
This means, that for $N \gg \lambda$, the `left' part of the convex domain dominates area and volume.
We also observe that
\begin{align*}
 |\Omega_N| &= N^3 |\Omega_1| + \mathcal{O}(N^2) \\
 |\partial \Omega_N| &= N^3 |\partial \Omega_1| + \mathcal{O}(N^2),
 \end{align*}
 where the implicit constants depend on $\Omega_1$ and $\Omega_2$.
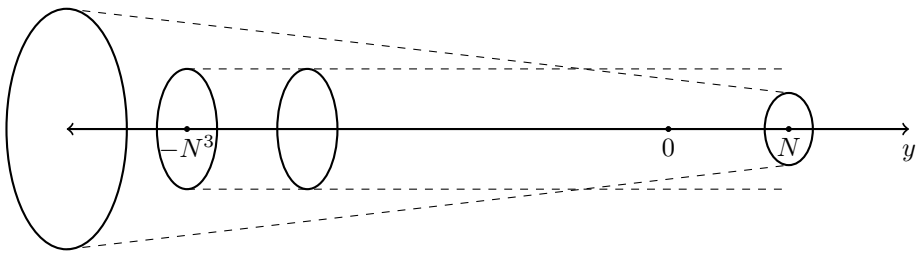
\begin{figure}[h!]
\begin{center}
\begin{tikzpicture}[scale=1.6]
\draw [<->, thick] (-1,0) -- (6,0);
\node at (6, -0.2) {$y$};
\draw[thick] (1,0) ellipse (0.25cm and 0.5cm);
\draw[thick] (0,0) ellipse (0.25cm and 0.5cm);
\filldraw (0,0) circle (0.02cm);
\node at (0, -0.15) {$-N^3$};
\draw [dashed] (0, 0.5) -- (5,0.5);
\draw [dashed] (0, -0.5) -- (5,-0.5);
\filldraw (4,0) circle (0.02cm);
\node at (4, -0.15) {0};
\draw[thick] (5,0) ellipse (0.2cm and 0.3cm);
\filldraw (5,0) circle (0.02cm);
\node at (5, -0.15) {$N$};
\draw[thick] (-1,0) ellipse (0.5cm and 1cm);
\draw[dashed] (-1, 1) -- (5, 0.3);
\draw[dashed] (-1, -1) -- (5, -0.3);
\end{tikzpicture}
\end{center}
\caption{The construction of $C_1$ and $C_2$.}
\end{figure}

Since $\Omega_2 \subset \Omega_1$, we have that
 $$ \Omega_N \cap \left\{(x,y) \in \mathbb{R}^{n}: y \geq 0 \right\} =  C_2 \cap \left\{(x,y) \in \mathbb{R}^{n}: y \geq 0 \right\}.$$
We now define a convex function on $\mathbb{R}^{n}$ via
$$ f(x,y) = \begin{cases} y \quad &\mbox{if}~y \geq 0, \\ 0 \quad &\mbox{otherwise.} \end{cases}$$
We obtain
$$ \int_{\Omega_N} f~ dx dy =  \int_{\Omega_N \cap \left\{y > 0\right\}} f~ dx dy=  \int_{C_2 \cap \left\{y > 0\right\}} f ~dx dy = (1+o(1)) \frac{N^2}{2} |\Omega_2|$$
$$ \int_{\partial \Omega_N} f~ d\sigma =   \int_{ \partial \Omega_N \cap \left\{y > 0\right\}} f~ d\sigma  = \int_{ \partial C_2 \cap \left\{y > 0\right\}} f ~d\sigma = (1+o(1)) \frac{N^2}{2} |\partial \Omega_2|.$$
This shows that
$$ \frac{1}{|\Omega_{N}|} \int_{\Omega_{N}}{f_{} ~ dx} = \frac{(1+o(1))}{2N} \frac{|\Omega_2|}{|\Omega_1|}$$
and
$$ \frac{1}{|\partial \Omega_{N}|} \int_{\partial \Omega_{N}}{f_{} ~ d\sigma} = \frac{(1+o(1))}{2N} \frac{|\partial \Omega_2|}{|\partial \Omega_1|}$$
which implies the desired result for $N \rightarrow \infty$.
\end{proof}

\section{Proof of the Geometric Lemma}

\begin{proof}
By the inequality 
$$\frac{|\Omega|}{|\partial\Omega|}\leq \mbox{inrad}(\Omega)\leq n \frac{|\Omega|}{|\partial\Omega|},$$
 a proof of which can be found in~\cite{larson2}, the supremum is no larger than $n$ since
\begin{equation}\label{eq:upper bound}
	\frac{|\partial \Omega|}{|\Omega|} \frac{| \Omega'|}{|\partial \Omega'|}\leq \frac{n \cdot \mbox{inrad}(\Omega')}{\mbox{inrad}(\Omega)}\leq n.
\end{equation}
What remains is to prove that this upper bound is saturated. 
The underlying idea of our proof is a theorem of Kawohl and Lachand-Robert \cite{kawohl} characterizing the Cheeger set of a convex set $\Omega\subset \mathbb{R}^2$. Specifically, their theorem states that for a convex $\Omega\subset \mathbb{R}^2$ the Cheeger problem
\begin{equation*}
	h(\Omega)=\inf\biggl\{\frac{|\partial\Omega'|}{|\Omega'|}: \Omega'\subset \overline{ \Omega}\biggr\}
\end{equation*}
is solved by the set 
$$
	\Omega' = \{x\in \Omega: \exists y \in \Omega \mbox{ such that } x\in B_{1/h(\Omega)^{}}(y)\subset \Omega\},
$$ 
where $B_{r}(x_0)$ is a ball of radius $r$ centered at $x_0$. 
We recall our use of the notation
\begin{equation*}
	\Omega_t = \{x\in \Omega: d(x, \partial\Omega)>t\},
\end{equation*}
where $d(x, \partial \Omega)$ denotes the distance to the boundary
$$ d(x, \partial \Omega) = \inf_{y \in \partial \Omega}{ \|x-y\|},$$
and we can equivalently write the Cheeger set of $\Omega$ as $\Omega'= \Omega_{1/h(\Omega)}+B_{1/h(\Omega)}$, where $B_r$ denotes a ball of radius $r$ centered in 0. Here and in what follows the sum of two sets is to be interpreted in the sense of the Minkowski sum:
\begin{equation*}
	A+B = \{x: x=a+b, a\in A, b\in B\}.
\end{equation*}

The situation when $n\geq 2$ is more complicated, and as far as we know a precise solution of the Cheeger problem is not available~\cite{alter}. Nonetheless, our aim in what follows is to prove that by taking $\Omega$ as a very thin $n$-simplex we can find a good enough candidate for $\Omega'$ among the one-parameter family of sets
\begin{equation}\label{eq:Kawohl_trialsest}
	   	\Omega_t + B_t, \quad 0\leq t\leq \mbox{inrad}(\Omega).
\end{equation}   

We construct our candidate for $\Omega$ as follows. Let $\Omega(\eta)\subset \mathbb{R}^n$ be the $n$-simplex obtained by taking a regular $(n-1)$-simplex of sidelength $\eta\gg 1$ in the hyperplane $\{x \in \mathbb{R}^n:x_1=-1\}$ with $(-1, 0, \ldots, 0)$ as center of mass and adding the last vertex at $(h(\eta), 0, \ldots, 0)$, where $h(\eta)$ is chosen so that $\mbox{inrad}(\Omega(\eta))=1$. Note that, as $\eta$ becomes large, $h(\eta)$ is approximately $1$ and $|\Omega(\eta)|\sim \eta^{d-1}$. \\

By construction $B_1\subset \Omega(\eta)$, and it is the unique unit ball of maximal radius contained in $\Omega(\eta)$. Moreover, the set $\Omega(\eta)$ is a tangential body to this ball (that is, a convex body all of whose supporting hyperplanes are tangential to the same ball). Since every tangential body to a ball is homothetic to its form body~\cite{Schneider1} (in our case $\Omega(\eta)$ is in fact equal to its form body), the main result in~\cite{larson1} implies 
\begin{equation*}
	|(\partial \Omega(\eta))_t| = (1-t)^{n-1}|\partial\Omega(\eta)|, \quad \mbox{ for all }t \in [0, 1].
\end{equation*}
An application of the coarea formula now yields the identity
\begin{equation}\label{eq:vol_per_ratio}
	|\Omega(\eta)|= \int_0^{1}|\partial(\Omega(\eta))_t|dt = \frac{|\partial\Omega(\eta)|}{n}.
\end{equation}
We also note that  $\Omega(\eta)\subset B_{2\eta}$. To see why this is true, we first note that the inradius of the regular $n$-simplex (by which we mean $n+1$ points all at distance 1 from each other embedded in $\mathbb{R}^n$) is given by
$$ r_n = \frac{1}{\sqrt{2n(n+1)}}.$$
The regular simplex is the convex body for which John's theorem is sharp, the circumradius is thus given by
$$ R_n = n\cdot r_n = \frac{\sqrt{n}}{\sqrt{2(n+1)}} \leq \frac{1}{\sqrt{2}}.$$
This shows that $\Omega(\eta) \subset B_{2\eta}$ (for the purpose of the proof, the constant 2 is not important and could be replaced by a much larger (absolute) constant).
Since it makes the computations somewhat simpler we consider, for a suitably chosen number $t$, the set $(1+t)\Omega(\eta)$. By construction $B_1\subset \Omega(\eta)$ which implies the inclusion 
$$
	\Omega(\eta)+B_t\subset \Omega(\eta)+t\Omega(\eta)= (1+t)\Omega(\eta).
$$
In particular, we can test~\eqref{convex} with $\Omega=(1+t)\Omega(\eta)$ and $\Omega'=\Omega(\eta)+B_t$ for any values of $t, \eta\gg 1$.
We note that up to rescaling by $(1+t)^{-1}$ this is exactly the family of sets in~\eqref{eq:Kawohl_trialsest}. Indeed, for each $t$ the set $((1+t)\Omega(\eta))_t=\Omega(\eta)$.\\

The final step of the proof is to show that by letting $t, \eta \to \infty$ appropriately 
\begin{equation}\label{eq:goal}
	\frac{|\partial((1+t)\Omega(\eta))|}{|(1+t)\Omega(\eta)|}\frac{|\Omega(\eta)+B_t|}{|\partial(\Omega(\eta)+B_t)|} \to n.
\end{equation}

To prove~\eqref{eq:goal} we recall the definition and some basic properties of mixed volumes~\cite[p. 275ff]{Schneider1}. Let $\mathcal{K}$ denote the set of convex bodies in $\mathbb{R}^n$ with nonempty interior. The mixed volume is defined as the unique symmetric function $W\colon \mathcal{K}^n \to \mathbb{R}_+$ satisfying
\begin{equation*}
  |\eta_1\Omega_1 + \ldots + \eta_m\Omega_m|= \sum_{j_1=1}^m\cdots \sum_{j_n=1}^m \eta_{j_1}\!\cdots \eta_{j_n}W(\Omega_{j_1}, \ldots, \Omega_{j_n}),
\end{equation*}
for any $\Omega_{1}, \ldots, \Omega_m \in \mathcal{K}$ and $\eta_1, \ldots, \eta_m\geq 0$~\cite{Schneider1}. Then $W$ satisfies the following properties:
\begin{enumerate}
  \item\label{itma} $W(\Omega_1, \ldots, \Omega_n)>0$ for $\Omega_1, \ldots, \Omega_n \in \mathcal{K}$.
  \item\label{itmb} $W$ is a multilinear function with respect to Minkowski addition.
  \item\label{itmc} $W$ is increasing with respect to inclusions in each of its arguments.
  \item\label{itme} The volume and perimeter of $\Omega\in \mathcal{K}$ can be written in terms of $W$:
  $$
    |\Omega|=W( \underbrace{\Omega, \ldots, \Omega}_{n~\mbox{\tiny times}}) \quad \mbox{and} \quad  |\partial\Omega|= n W(\underbrace{\Omega, \ldots, \Omega}_{n-1~\mbox{\tiny times}}, B_1).
  $$
\end{enumerate}
Since only mixed volumes of two distinct sets appear in our proof we introduce the shorthand notation
$$
	W_j(K, L) = W(\underbrace{K, \ldots, K}_{n-j}, \underbrace{L, \ldots, L}_{j}).
$$
By~\eqref{eq:vol_per_ratio}, we have
$$ 	\frac{|\partial((1+t)\Omega(\eta))|}{|(1+t)\Omega(\eta)|} = \frac{n}{1+t}.$$
The definition of $W$ implies
$$ |\Omega(\eta)+B_t| = W(\Omega(\eta)+B_t, \dots, \Omega(\eta)+B_t).$$
Multilinearity allows us to expand this term as 
$$ W(\Omega(\eta)+B_t, \dots, \Omega(\eta)+B_t) = \sum_{j=0}^n \binom{n}{j}t^j W_j(\Omega(\eta), B_1)$$
and the same argument shows
$$ |\partial(\Omega(\eta)+B_t)| = n\sum_{j=0}^{n-1}\binom{n-1}{j}t^jW_{j+1}(\Omega(\eta), B_1).$$
Altogether, we can write the expression of interest as
\begin{equation}\label{eq:expansion in W}
	\frac{|\partial((1+t)\Omega(\eta))|}{|(1+t)\Omega(\eta)|}\frac{|\Omega(\eta)+B_t|}{|\partial(\Omega(\eta)+B_t)|}
	=
	\frac{n}{(1+t)}\frac{\sum_{j=0}^n \binom{n}{j}t^j W_j(\Omega(\eta), B_1)}{n\sum_{j=0}^{n-1}\binom{n-1}{j}t^jW_{j+1}(\Omega(\eta), B_1)}.
\end{equation}
In order to prove~\eqref{eq:goal} we need a bound from below. For the sum in the numerator it suffices  to keep the first two terms in the expansion and to use~Property (\ref{itma}) of $W$ resulting in
\begin{align*}
 	\sum_{j=0}^n \binom{n}{j}t^j W_j(\Omega(\eta), B_1)&\geq W_0(\Omega(\eta), B_1) + n t W_1(\Omega(\eta), B_1) \\
 	&= |\Omega(\eta)|+ t|\partial\Omega(\eta)|.
\end{align*}
To bound the sum in the denominator we wish to keep the term with $j=0$ as is. 
For $j\geq 1$, we now use that $\Omega(\eta)\subset B_{2\eta}$ together with Property (\ref{itmc}) to bound
$$
	W_{j+1}(\Omega(\eta), B_1) \leq W_{j+1}(B_{2\eta}, B_{1})=(2\eta)^{n-j-1}|B_1|.
$$
Inserting the two bounds above into~\eqref{eq:expansion in W} yields
\begin{align*}
\frac{|\partial((1+t)\Omega(\eta))|}{|(1+t)\Omega(\eta)|}&\frac{|\Omega(\eta)+B_t|}{|\partial(\Omega(\eta)+B_t)|}\\
	&\geq
	\frac{n}{(1+t)}\frac{|\Omega(\eta)|+t|\partial\Omega(\eta)|}{|\partial\Omega(\eta)|+n\sum_{j=1}^{n-1}\binom{n-1}{j}2^{n-j-1}t^j\eta^{n-j-1} |B_1|}.
\end{align*}
We recall that $|\Omega(\eta)| = n^{-1} |\partial \Omega(\eta)|$ and therefore
\begin{align*}
\frac{|\partial((1+t)\Omega(\eta))|}{|(1+t)\Omega(\eta)|}&\frac{|\Omega(\eta)+B_t|}{|\partial(\Omega(\eta)+B_t)|}\\
	&\geq
	\frac{n t}{(1+t)}\frac{1 + 1/(nt)}{1+n|B_1|\sum_{j=1}^{n-1}\binom{n-1}{j}2^{n-j-1}t^j\eta^{n-j-1}/ |\partial\Omega(\eta)| }.
\end{align*}

By construction $|\partial\Omega(\eta)|\sim \eta^{n-1}$ as $\eta \to \infty$. Consequently, by choosing $t=\sqrt{\eta}$ (though the more general choice $t=\eta^{\alpha}$ for $0<\alpha<1$ would also work) we find $t^j \eta^{n-j-1}/|\partial\Omega(\eta)|\sim \eta^{-j/2}$. Therefore, taking $\eta$ (and thus $t$) to infinity, we obtain
\begin{equation*}
	\liminf_{\eta \to \infty}\frac{|\partial((1+\sqrt{\eta})\Omega(\eta))|}{|(1+\sqrt{\eta})\Omega(\eta)|}\frac{|\Omega(\eta)+B_{\sqrt{\eta}}|}{|\partial(\Omega(\eta)+B_{\sqrt{\eta}})|} \geq n,
\end{equation*}
which when combined with the matching upper bound~\eqref{eq:upper bound} completes the proof.
\end{proof}

\section{Proof of Theorem 3}
\begin{proof}
We use the estimate
$$ \int_{\Omega}{ f dx} \leq \max_{x \in \partial \Omega}{  \frac{\partial u}{\partial \nu}(x)} \int_{\partial \Omega} { f d\sigma}$$
introduced in Proposition 1 above and use, inspired by the argument in \cite{jianfeng}, estimates for the torsion function. One such estimate for the torsion function comes from $P-$functions, we refer to the classic book of Sperb \cite[Eq. (6.12)]{sperb},
$$  \max_{x \in \partial \Omega}{  \frac{\partial u}{\partial \nu}(x)} \leq \sqrt{2} \|u\|^{\frac{1}{2}}_{L^{\infty}}.$$
It remains to estimate the largest value of the torsion function. There are two different approaches: we can interpret it as the maximum lifetime of Brownian motion inside a domain of given measure or we can interpret it as the solution of a partial differential equation to which Talenti's theorem \cite{talenti} can be applied. In both cases, we end up with a standard isoperimetric estimate \cite{talenti} (that was also used in \cite{jianfeng})
$$ \|u\|_{L^{\infty}} \leq \frac{1}{2n}\left( \frac{|\Omega|}{\omega_n}\right)^{\frac{2}{n}}$$
to obtain
$$  \max_{x \in \partial \Omega}{  \frac{\partial u}{\partial \nu}(x)}  \leq \frac{ |\Omega|^{1/n}}{\omega_n^{1/n} \sqrt{n}}.$$
\end{proof}

\textbf{Acknowledgment.} This research was initiated at the workshop `Shape Optimization with Surface Interactions' at the American Institute of Mathematics in June 2019. The authors are grateful 
to the organizers of the workshop as well as the Institute. KB's research was supported in part by Simons Foundation Grant 506732. 
JL’s research was supported in part by a Bucknell University Scholarly Development Grant.
SL acknowledges financial support from Swedish Research Council Grant No. 2012-3864. SS's research was supported in part by the NSF (DMS-1763179) and the Alfred P. Sloan foundation.

\end{document}